\documentclass{amsart}
\usepackage{amssymb}
\usepackage{amsfonts}

\newtheorem{theorem}{Theorem}[section]
\theoremstyle{plain}

\newtheorem{corollary}{Corollary}[section]

\newtheorem{definition}{Definition}[section]
\newtheorem{example}{Example}[section]

\newtheorem{lemma}{Lemma}[section]

\newtheorem{proposition}{Proposition}[section]

\numberwithin{equation}{section}
\input{tcilatex}

\begin{document}
\title[Pairwise Semiregular Properties]{Pairwise Semiregular Properties on
Generalized Pairwise Lindel\"{o}f Spaces}
\author{Zabidin Salleh}
\address{School of Informatics and Applied Mathematics, Universiti Malaysia
Terengganu, 21030 Kuala Nerus, Terengganu, Malaysia.}
\email{zabidin@umt.edu.my}
\date{January 28, 2019}

\begin{abstract}
Let $\left( X,\tau _{1},\tau _{2}\right) $ be a bitopological space and $%
\left( X,\tau _{\left( 1,2\right) }^{s},\tau _{\left( 2,1\right)
}^{s}\right) $ its pairwise semiregularization. Then a bitopological
property $\mathcal{P}$\ is called pairwise semiregular provided that $\left(
X,\tau _{1},\tau _{2}\right) $\ has the property $\mathcal{P}$\ if and only
if $\left( X,\tau _{\left( 1,2\right) }^{s},\tau _{\left( 2,1\right)
}^{s}\right) $\ has the same property. In this work we study pairwise
semiregular property of  $\left( i,j\right) $-nearly Lindel\"{o}f, pairwise
nearly Lindel\"{o}f, $\left( i,j\right) $-almost Lindel\"{o}f, pairwise
almost Lindel\"{o}f, $\left( i,j\right) $-weakly Lindel\"{o}f and pairwise
weakly Lindel\"{o}f spaces. We prove that $\left( i,j\right) $-almost Lindel%
\"{o}f, pairwise almost Lindel\"{o}f, $\left( i,j\right) $-weakly Lindel\"{o}%
f and pairwise weakly Lindel\"{o}f are pairwise semiregular properties, on
the contrary of each type of pairwise Lindel\"{o}f space which are not
pairwise semiregular properties.
\end{abstract}

\maketitle

\section{\textbf{Introduction}}

\noindent Semiregular properties in topological spaces have been studied by
many topologist. Some of them related to this research studied by Mrsevic et
al. \cite{MrsevicReillyVamanamurthy85, MrsevicReillyVamanamurthy86} and
Fawakhreh and K\i l\i \c{c}man \cite{FawaAdem04}. The purpose of this paper
is to study pairwise semiregular properties on generalized pairwise Lindel%
\"{o}f spaces, that we have studied in \cite{BidiAdem, AdemBidi,
BidiAdem2014, BidiAdem1}, namely, $\left( i,j\right) $-nearly Lindel\"{o}f,
pairwise nearly Lindel\"{o}f, $\left( i,j\right) $-almost Lindel\"{o}f,
pairwise almost Lindel\"{o}f, $\left( i,j\right) $-weakly Lindel\"{o}f and
pairwise weakly Lindel\"{o}f spaces.

The main results is that the Lindel\"{o}f, $B$-Lindel\"{o}f, $s$-Lindel\"{o}%
f and $p$-Lindel\"{o}f spaces are not pairwise semiregular properties. While 
$\left( i,j\right) $-almost Lindel\"{o}f, pairwise almost Lindel\"{o}f, $%
\left( i,j\right) $-weakly Lindel\"{o}f and pairwise weakly Lindel\"{o}f
spaces are pairwise semiregular properties. We also show that $\left(
i,j\right) $-nearly Lindel\"{o}f and pairwise nearly Lindel\"{o}f spaces are
satisfying pairwise semiregular invariant properties.

\section{\textbf{Preliminaries}}

\noindent Throughout this paper, all spaces $\left( X,\tau \right) $ and $%
\left( X,\tau _{1},\tau _{2}\right) $ $($or simply $X)$ are always mean
topological spaces and bitopological spaces, respectively unless explicitly
stated. If $\mathcal{P}$ is a topological property, then $\left( \tau
_{i},\tau _{j}\right) $-$\mathcal{P}$ denotes an analogue of this property
for $\tau _{i}$ has property $\mathcal{P}$ with respect to $\tau _{j}$, and $%
p$-$\mathcal{P}$ denotes the conjunction $\left( \tau _{1},\tau _{2}\right) $%
-$\mathcal{P}\wedge \left( \tau _{2},\tau _{1}\right) $-$\mathcal{P}$, i.e., 
$p$-$\mathcal{P}$ denotes an absolute bitopological analogue of $\mathcal{P}$%
. As we shall see below, sometimes $\left( \tau _{1},\tau _{2}\right) $-$%
\mathcal{P}\Longleftrightarrow \left( \tau _{2},\tau _{1}\right) $-$\mathcal{%
P}$ (and thus $\Leftrightarrow p$-$\mathcal{P}$) so that it suffices to
consider one of these three bitopological analogue. Also sometimes $\tau
_{1} $-$\mathcal{P}\Longleftrightarrow \tau _{2}$-$\mathcal{P}$ and thus $%
\mathcal{P}\Longleftrightarrow \tau _{1}$-$\mathcal{P}\wedge \tau _{2}$-$%
\mathcal{P}$, i.e., $\left( X,\tau _{i}\right) $ has property $\mathcal{P}$
for each $i=1,2$.

Also note that $\left( X,\tau _{i}\right) $ has a property $\mathcal{P}%
\Longleftrightarrow \left( X,\tau _{1},\tau _{2}\right) $ has a property $%
\tau _{i}$-$\mathcal{P}$. Sometimes the prefixes $\left( \tau _{i},\tau
_{j}\right) $- or $\tau _{i}$- will be replaced by $\left( i,j\right) $- or $%
i$- respectively, if there is no chance for confusion. By $i$-open cover of $%
X$, we mean that the cover of $X$ by $i$-open sets in $X$; similar for the $%
\left( i,j\right) $-regular open cover of $X$ etc. By $i$-$\limfunc{int}%
\left( A\right) $ and $i$-$\limfunc{cl}\left( A\right) $, we shall mean the
interior and the closure of a subset $A$ of $X$ with respect to topology $%
\tau _{i}$, respectively. In this paper always $i,j\in \left\{ 1,2\right\} $
and $i\neq j$. The reader may consult \cite{Dvalish05} for the detail
notations.

The following are some basic concepts.

\begin{definition}
\cite{KhedShib91, SingAry71}\textbf{\ }A subset $S$\ of a bitopological
space $\left( X,\tau _{1},\tau _{2}\right) $\ is said to be $\left(
i,j\right) $-regular open $($resp. $\left( i,j\right) $-regular closed$)$\
if $i$-$\limfunc{int}\left( j\text{-}\limfunc{cl}\left( S\right) \right) =S$%
\ $($resp. $i$-$\limfunc{cl}\left( j\text{-}\limfunc{int}\left( S\right)
\right) =S)$, where $i,j\in \left\{ 1,2\right\} ,i\neq j$. $S$\ is called
pairwise regular open $($resp. pairwise regular closed$)$\ if it is both $%
\left( 1,2\right) $-regular open and $\left( 2,1\right) $-regular open $($%
resp. $\left( 1,2\right) $-regular closed and $\left( 2,1\right) $-regular
closed$)$.
\end{definition}

\begin{definition}
\cite{KhedShib91, SingSing70}\textbf{\ }A bitopological space $\left( X,\tau
_{1},\tau _{2}\right) $\ is said to be $\left( i,j\right) $-almost regular
if for each point $x\in X$\ and for each $\left( i,j\right) $-regular open\
set $V$\ containing $x$, there exists an $\left( i,j\right) $-regular open
set $U$\ such that $x\in U\subseteq j$-$\limfunc{cl}\left( U\right)
\subseteq V$. $X$\ is called pairwise almost regular if it is both $\left(
1,2\right) $-almost regular\ and $\left( 2,1\right) $-almost regular.
\end{definition}

In any bitopological space $\left( X,\tau _{1},\tau _{2}\right) $, the
family of all $\left( i,j\right) $-regular open sets is closed under finite
intersections. Thus the family of $\left( i,j\right) $-regular open sets in
any bitopological space $\left( X,\tau _{1},\tau _{2}\right) $ forms a base
for a coarser topology called $\left( i,j\right) $-semiregularization of $%
\left( X,\tau _{1},\tau _{2}\right) $, which is defined as follows.\ 

\begin{definition}
\cite{BidiAdem} The topology generated by the $\left( i,j\right) $-regular
open subsets of $\left( X,\tau _{1},\tau _{2}\right) $\ is denoted by $\tau
_{\left( i,j\right) }^{s}$\ and it is called $\left( i,j\right) $%
-semiregularization of $X$. The topologies is pairwise semiregularization of 
$X$ if the first topology is $\left( 1,2\right) $-semiregularization of $X$\
and the second topology is $\left( 2,1\right) $-semiregularization of $X$.
If $\tau _{i}\equiv \tau _{\left( i,j\right) }^{s}$, then $X$\ is said to be 
$\left( i,j\right) $-semiregular. $\left( X,\tau _{1},\tau _{2}\right) $\ is
called pairwise semiregular if it is both $\left( 1,2\right) $-semiregular
and $\left( 2,1\right) $-semiregular, that is, whenever $\tau _{i}\equiv
\tau _{\left( i,j\right) }^{s}$\ for each $i,j\in \left\{ 1,2\right\} $\ and 
$i\neq j$. In other words, $\left( X,\tau _{1},\tau _{2}\right) $ is $\left(
i,j\right) $-semiregular if the family of $\left( i,j\right) $-regular open
sets form a base for the topology $\tau _{i}$.
\end{definition}

It is very clear that $\tau _{\left( i,j\right) }^{s}\subseteq \tau _{i}$,
but it is not necessary $\tau _{i}\subseteq \tau _{\left( i,j\right) }^{s}$.
Thus with every given bitopological space $\left( X,\tau _{1},\tau
_{2}\right) $ there is associated another bitopological space $\left( X,\tau
_{\left( 1,2\right) }^{s},\tau _{\left( 2,1\right) }^{s}\right) $ in the
manner described above (see \cite{SingAry71}). We provide the following
example in order to understand the concept of pairwise semiregular spaces
clearly.

\begin{example}
For the set of all real numbers $%
\mathbb{R}
$, let $\tau _{u}$\ denotes the usual topology and $\tau _{s}$\ denote the
Sorgenfrey topology, i.e., topology generated by right half-open intervals $%
( $see \cite{Steen78}$)$. Then $\left( 
\mathbb{R}
,\tau _{u},\tau _{s}\right) $\ is $\left( \tau _{u},\tau _{s}\right) $%
-semiregular since $\tau _{u}=\tau _{\left( \tau _{u},\tau _{s}\right) }^{s}$%
, i.e., $\tau _{u}$\ generated by $\left( \tau _{u},\tau _{s}\right) $%
-regular open subsets of $%
\mathbb{R}
$. $\left( 
\mathbb{R}
,\tau _{u},\tau _{s}\right) $\ is also $\left( \tau _{s},\tau _{u}\right) $%
-semiregular since $\tau _{s}=\tau _{\left( \tau _{s},\tau _{u}\right) }^{s}$%
\ because any set $E\in \tau _{s}$ is the union of a collection of $\left(
\tau _{s},\tau _{u}\right) $-regular open sets in $%
\mathbb{R}
$. Thus $\left( 
\mathbb{R}
,\tau _{u},\tau _{s}\right) $\ is pairwise semiregular.
\end{example}

Khedr and Alshibani \cite{KhedShib91} defined the equivalent definition of $%
\left( i,j\right) $-semiregular spaces as follows.

\begin{definition}
A bitopological space $X$\ is said to be $\left( i,j\right) $-semiregular if
for each $x\in X$\ and for each $i$-open subset $V$\ of $X$\ containing $x$,
there is an $i$-open set $U$\ such that $x\in U\subseteq i$-$\limfunc{int}%
\left( j\text{-}\limfunc{cl}\left( U\right) \right) \subseteq V$. $X$\ is
called pairwise semiregular if it is both $\left( 1,2\right) $-semiregular
and $\left( 2,1\right) $-semiregular.
\end{definition}

\begin{definition}
\cite{Datta72} A bitopological space $\left( X,\tau _{1},\tau _{2}\right) $\
is said to be $\left( i,j\right) $-extremally disconnected if the $i$%
-closure of every $j$-open set is $j$-open. $X$\ is called pairwise
extremally disconnected if it is both $\left( 1,2\right) $-extremally
disconnected and $\left( 2,1\right) $-extremally disconnected.
\end{definition}

Recall that a property $\mathcal{P}$ will be called bitopological property
(resp. $p$-topological property) if whenever $\left( X,\tau _{1},\tau
_{2}\right) $ has property $\mathcal{P}$, then every space homeomorphic
(resp. $p$-homeomorphic) to $\left( X,\tau _{1},\tau _{2}\right) $ also has
property $\mathcal{P}$ (see \cite{AdemBidi07-1}). If a bitopological space $%
X $ has bitopological (or $p$-topological) property $\mathcal{P}$, one may
ask, does the pairwise semiregularization of $X$ satisfies the property $%
\mathcal{P}$ also? Now we arrive to the concept of pairwise semiregular
property.\ 

\begin{definition}
\label{Def 2.6}Let $\left( X,\tau _{1},\tau _{2}\right) $\ be a
bitopological space and let $\left( X,\tau _{\left( 1,2\right) }^{s},\tau
_{\left( 2,1\right) }^{s}\right) $\ its pairwise semiregularization. A
bitopological property $\mathcal{P}$\ is called pairwise semiregular
provided that $\left( X,\tau _{1},\tau _{2}\right) $\ has the property $%
\mathcal{P}$\ if and only if $\left( X,\tau _{\left( 1,2\right) }^{s},\tau
_{\left( 2,1\right) }^{s}\right) $\ has the property $\mathcal{P}$.
\end{definition}

\begin{lemma}
\label{Lem 2.1}\cite{BidiAdem} Let $\left( X,\tau _{1},\tau _{2}\right) $\
be a bitopological space and let $\left( X,\tau _{\left( 1,2\right)
}^{s},\tau _{\left( 2,1\right) }^{s}\right) $\ its pairwise
semiregularization. Then
\end{lemma}

\noindent $\left( a\right) $\ $\tau _{i}$-$\limfunc{int}\left( C\right)
=\tau _{\left( i,j\right) }^{s}$-$\limfunc{int}\left( C\right) $\ for every $%
\tau _{j}$-closed set $C;$

\noindent $\left( b\right) $\ $\tau _{i}$-$\limfunc{cl}\left( A\right) =\tau
_{\left( i,j\right) }^{s}$-$\limfunc{cl}\left( A\right) $\ for every $A\in
\tau _{j};$

\noindent $\left( c\right) $\ the family of $\left( \tau _{i},\tau
_{j}\right) $-regular open sets of $\left( X,\tau _{1},\tau _{2}\right) $\
are the same as the family of $\left( \tau _{\left( i,j\right) }^{s},\tau
_{\left( j,i\right) }^{s}\right) $-regular open sets of $\left( X,\tau
_{\left( 1,2\right) }^{s},\tau _{\left( 2,1\right) }^{s}\right) ;$

\noindent $\left( d\right) $\ the family of $\left( \tau _{i},\tau
_{j}\right) $-regular closed sets of $\left( X,\tau _{1},\tau _{2}\right) $\
are the same as the family of $\left( \tau _{\left( i,j\right) }^{s},\tau
_{\left( j,i\right) }^{s}\right) $-regular closed sets of $\left( X,\tau
_{\left( 1,2\right) }^{s},\tau _{\left( 2,1\right) }^{s}\right) ;$

\noindent $\left( e\right) $ $\left( \tau _{\left( i,j\right) }^{s}\right)
_{\left( i,j\right) }^{s}=\tau _{\left( i,j\right) }^{s}$.

\section{\textbf{Pairwise Semiregularization of Pairwise Lindel\"{o}f Spaces 
\label{Sect 3}}}

\begin{definition}
\cite{ForaHdeib83, AdemBidi07}. A bitopological space $\left( X,\tau
_{1},\tau _{2}\right) $\ is said to be $i$-Lindel\"{o}f if the topological
space $\left( X,\tau _{i}\right) $\ is Lindel\"{o}f. $X$\ is called Lindel%
\"{o}f if it is $i$-Lindel\"{o}f for each $i=1,2$. In other words, $\left(
X,\tau _{1},\tau _{2}\right) $\ is called Lindel\"{o}f if the topological
space $\left( X,\tau _{1}\right) $\ and $\left( X,\tau _{2}\right) $\ are
both Lindel\"{o}f.
\end{definition}

Note that $i$-Lindel\"{o}f property as well as Lindel\"{o}f property is not
a pairwise semiregular property by the following example.\ 

\begin{example}
\label{Ex 3.1}Let $X$\ be a set with cardinality $2^{c}$, where $c=\limfunc{%
card}\left( 
\mathbb{R}
\right) $. Let $\tau _{1}$\ be a co-$c$\ topology on $X$\ consisting of $%
\emptyset $\ and all subsets of $X$\ whose complements have cardinality at
most $c$\ and let $\tau _{2}$\ be a cofinite topology on $X$. Then $\left(
X,\tau _{1},\tau _{2}\right) $\ is $\tau _{2}$-Lindel\"{o}f but not $\tau
_{1}$-Lindel\"{o}f and hence not Lindel\"{o}f. Observe that $\left( X,\tau
_{\left( 1,2\right) }^{s},\tau _{\left( 2,1\right) }^{s}\right) $\ is $\tau
_{\left( 1,2\right) }^{s}$-Lindel\"{o}f and $\tau _{\left( 2,1\right) }^{s}$%
-Lindel\"{o}f since $\tau _{\left( 1,2\right) }^{s}$\ and $\tau _{\left(
2,1\right) }^{s}$\ are indiscrete topologies. Hence $\left( X,\tau _{\left(
1,2\right) }^{s},\tau _{\left( 2,1\right) }^{s}\right) $\ is Lindel\"{o}f%
\textit{.}
\end{example}

\begin{definition}
A bitopological space $\left( X,\tau _{1},\tau _{2}\right) $\ is called $%
\left( i,j\right) $-Lindel\"{o}f \cite{ForaHdeib83, AdemBidi07}\ if for
every $i$-open cover of $X$\ there is a countable\ $j$-open subcover. $X$\
is called $B$-Lindel\"{o}f \cite{ForaHdeib83}\ or $p_{1}$-Lindel\"{o}f \cite%
{AdemBidi07}\ if it is both $\left( 1,2\right) $-Lindel\"{o}f\ and $\left(
2,1\right) $-Lindel\"{o}f.
\end{definition}

An $\left( i,j\right) $-Lindel\"{o}f property as well as $B$-Lindel\"{o}f
property is not pairwise semiregular property by the following example.

\begin{example}
Let $\left( X,\tau _{1},\tau _{2}\right) $\ be a bitopological space as in
Example \ref{Ex 3.1}. Then $\left( X,\tau _{1},\tau _{2}\right) $\ is not $%
\left( \tau _{1},\tau _{2}\right) $-Lindel\"{o}f but it is $\left( \tau
_{2},\tau _{1}\right) $-Lindel\"{o}f and hence not $B$-Lindel\"{o}f. Observe
that $\left( X,\tau _{\left( 1,2\right) }^{s},\tau _{\left( 2,1\right)
}^{s}\right) $\ is $\left( \tau _{\left( 1,2\right) }^{s},\tau _{\left(
2,1\right) }^{s}\right) $-Lindel\"{o}f and $\left( \tau _{\left( 2,1\right)
}^{s},\tau _{\left( 1,2\right) }^{s}\right) $-Lindel\"{o}f since $\tau
_{\left( 1,2\right) }^{s}$\ and $\tau _{\left( 2,1\right) }^{s}$\ are
indiscrete topologies. Hence $\left( X,\tau _{\left( 1,2\right) }^{s},\tau
_{\left( 2,1\right) }^{s}\right) $\ is $B$-Lindel\"{o}f.
\end{example}

\begin{definition}
A cover $U\ $of a bitopological space $\left( X,\tau _{1},\tau _{2}\right) $%
\ is called $\tau _{1}\tau _{2}$-open \cite{Swart71} if $U\subseteq \tau
_{1}\cup \tau _{2}$. If, in addition, $U$\ contains at least one nonempty
member of $\tau _{1}$\ and at least one nonempty member of $\tau _{2}$, it
is called $p$-open \cite{Fletcher69}.
\end{definition}

\begin{definition}
\cite{ForaHdeib83}\textbf{\ }A bitopological space $\left( X,\tau _{1},\tau
_{2}\right) $\ is called $s$-Lindel\"{o}f \ $($resp. $p$-Lindel\"{o}f$)$\ if
every $\tau _{1}\tau _{2}$-open $($resp. $p$-open$)$\ cover of $X$\ has a
countable\ subcover.
\end{definition}

A $p$-Lindel\"{o}f property is not pairwise semiregular property by the
following example. Thus the $s$-Lindel\"{o}f property is also not pairwise
semiregular property.

\begin{example}
Let $\left( X,\tau _{1},\tau _{2}\right) $\ be a bitopological space as in
Example \ref{Ex 3.1}. Then $\left( X,\tau _{1},\tau _{2}\right) $\ is not $p$%
-Lindel\"{o}f and hence not $s$-Lindel\"{o}f. Observe that $\left( X,\tau
_{\left( 1,2\right) }^{s},\tau _{\left( 2,1\right) }^{s}\right) $\ is $p$%
-Lindel\"{o}f and $s$-Lindel\"{o}f since $\tau _{\left( 1,2\right) }^{s}$\
and $\tau _{\left( 2,1\right) }^{s}$\ are indiscrete topologies.
\end{example}

\section{\textbf{Pairwise Semiregularization of Generalized Pairwise Lindel%
\"{o}f Spaces}}

\begin{definition}
\cite{BidiAdem, AdemBidi, BidiAdem1}\textbf{\ }A bitopological space $X$\ is
said to be $\left( i,j\right) $-nearly Lindel\"{o}f $($resp. $\left(
i,j\right) $-almost Lindel\"{o}f, $\left( i,j\right) $-weakly Lindel\"{o}f$)$%
\ if for every $i$-open cover $\left\{ U_{\alpha }:\alpha \in \Delta
\right\} $\ of $X$, there exists a countable subset $\left\{ \alpha
_{n}:n\in 
\mathbb{N}
\right\} $\ of $\Delta $\ such that $X=\dbigcup\limits_{n\in 
\mathbb{N}
}i$-$\limfunc{int}\left( j\text{-}\limfunc{cl}\left( U_{\alpha _{n}}\right)
\right) $\ $\left( \text{\textit{resp.}}\mathit{\ }X=\dbigcup\limits_{n\in 
\mathbb{N}
}j\text{\textit{-}}\limfunc{cl}\left( U_{\alpha _{n}}\right) \text{\textit{, 
}}X=j\text{\textit{-}}\limfunc{cl}\left( \dbigcup\limits_{n\in 
\mathbb{N}
}\left( U_{\alpha _{n}}\right) \right) \right) $.\ $X$\ is called pairwise
nearly Lindel\"{o}f $($resp. pairwise almost Lindel\"{o}f, pairwise weakly
Lindel\"{o}f$)$\ if it\ is both $\left( 1,2\right) $-nearly Lindel\"{o}f $($%
resp. $\left( 1,2\right) $-almost Lindel\"{o}f, $\left( 1,2\right) $-weakly
Lindel\"{o}f$)$\ and $\left( 2,1\right) $-nearly Lindel\"{o}f $($resp. $%
\left( 2,1\right) $-almost Lindel\"{o}f, $\left( 2,1\right) $-weakly Lindel%
\"{o}f$)$.
\end{definition}

Our first result is analogue with the result of Mr\v{s}evi\'{c} et al. \cite[%
Theorem 1]{MrsevicReillyVamanamurthy86}.

\begin{theorem}
\label{Th 4.1}A bitopological space $\left( X,\tau _{1},\tau _{2}\right) $
is $\left( \tau _{i},\tau _{j}\right) $-nearly Lindel\"{o}f if and only if $%
\left( X,\tau _{\left( 1,2\right) }^{s},\tau _{\left( 2,1\right)
}^{s}\right) $\ is $\tau _{\left( i,j\right) }^{s}$-Lindel\"{o}f.
\end{theorem}

\begin{proof}
Let $\left( X,\tau _{1},\tau _{2}\right) $ be a $\left( \tau _{i},\tau
_{j}\right) $-nearly Lindel\"{o}f and let $\left\{ U_{\alpha }:\alpha \in
\Delta \right\} $ be a $\tau _{\left( i,j\right) }^{s}$-open cover of $%
\left( X,\tau _{\left( 1,2\right) }^{s},\tau _{\left( 2,1\right)
}^{s}\right) $. For each $x\in X$, there exists $\alpha _{x}\in \Delta $
such that $x\in U_{\alpha _{x}}$ and since for each $\alpha _{x}\in \Delta
,U_{\alpha _{x}}\in $ $\tau _{\left( i,j\right) }^{s}$, there exists a $%
\left( \tau _{i},\tau _{j}\right) $-regular open set $V_{\alpha _{x}}$ in $%
\left( X,\tau _{1},\tau _{2}\right) $ such that $x\in V_{\alpha
_{x}}\subseteq U_{\alpha _{x}}$. So $X=\dbigcup\nolimits_{x\in X}V_{\alpha
_{x}}$ and hence $\left\{ V_{\alpha _{x}}:x\in X\right\} $ is a $\left( \tau
_{i},\tau _{j}\right) $-regular open cover of $X$. Since $\left( X,\tau
_{1},\tau _{2}\right) $ is $\left( \tau _{i},\tau _{j}\right) $-nearly Lindel%
\"{o}f, there exists a countable subset of points $x_{1},\ldots
,x_{n},\ldots $ of $X$ such that $X=\dbigcup\nolimits_{n\in 
\mathbb{N}
}V_{\alpha _{x_{n}}}\subseteq \dbigcup\nolimits_{n\in 
\mathbb{N}
}U_{\alpha _{x_{n}}}$. This shows that $\left( X,\tau _{\left( 1,2\right)
}^{s},\tau _{\left( 2,1\right) }^{s}\right) $ is $\tau _{\left( i,j\right)
}^{s}$-Lindel\"{o}f.

Conversely, suppose that $\left( X,\tau _{\left( 1,2\right) }^{s},\tau
_{\left( 2,1\right) }^{s}\right) $\ is $\tau _{\left( i,j\right) }^{s}$%
-Lindel\"{o}f and let $\left\{ V_{\alpha }:\alpha \in \Delta \right\} $ be a 
$\left( \tau _{i},\tau _{j}\right) $-regular open cover of $\left( X,\tau
_{1},\tau _{2}\right) $. Since $V_{\alpha }\in \tau _{\left( i,j\right)
}^{s} $ for each $\alpha \in \Delta ,\left\{ V_{\alpha }:\alpha \in \Delta
\right\} $ is a $\tau _{\left( i,j\right) }^{s}$-open cover of $\left(
X,\tau _{\left( 1,2\right) }^{s},\tau _{\left( 2,1\right) }^{s}\right) $.
Since $\left( X,\tau _{\left( 1,2\right) }^{s},\tau _{\left( 2,1\right)
}^{s}\right) $\ is $\tau _{\left( i,j\right) }^{s}$-Lindel\"{o}f, there
exists a countable subcover such that $X=\dbigcup\nolimits_{n\in 
\mathbb{N}
}V_{\alpha _{n}}$. This implies that $\left( X,\tau _{1},\tau _{2}\right) $
is $\left( \tau _{i},\tau _{j}\right) $-nearly Lindel\"{o}f.
\end{proof}

\begin{corollary}
A bitopological space $\left( X,\tau _{1},\tau _{2}\right) $ is pairwise
nearly Lindel\"{o}f if and only if $\left( X,\tau _{\left( 1,2\right)
}^{s},\tau _{\left( 2,1\right) }^{s}\right) $\ is Lindel\"{o}f.
\end{corollary}

\begin{proposition}
\label{Prop 4.1}A bitopological space $\left( X,\tau _{\left( 1,2\right)
}^{s},\tau _{\left( 2,1\right) }^{s}\right) $ is $\left( \tau _{\left(
i,j\right) }^{s},\tau _{\left( j,i\right) }^{s}\right) $-nearly Lindel\"{o}f
if and only if $\left( X,\tau _{\left( 1,2\right) }^{s},\tau _{\left(
2,1\right) }^{s}\right) $\ is $\tau _{\left( i,j\right) }^{s}$-Lindel\"{o}f.
\end{proposition}

\begin{proof}
The sufficient condition is obvious by the definitions. So we need only to
prove necessary condition. Suppose that $\left\{ U_{\alpha }:\alpha \in
\Delta \right\} $ is a $\tau _{\left( i,j\right) }^{s}$-open cover of $%
\left( X,\tau _{\left( 1,2\right) }^{s},\tau _{\left( 2,1\right)
}^{s}\right) $. For each $x\in X$, there exists $\alpha _{x}\in \Delta $
such that $x\in U_{\alpha _{x}}$. Since $\left( X,\tau _{\left( 1,2\right)
}^{s},\tau _{\left( 2,1\right) }^{s}\right) $ is $\left( \tau _{\left(
i,j\right) }^{s},\tau _{\left( j,i\right) }^{s}\right) $-semiregular, there
exists a $\tau _{\left( i,j\right) }^{s}$-open set $V_{\alpha _{x}}$ in $%
\left( X,\tau _{\left( 1,2\right) }^{s},\tau _{\left( 2,1\right)
}^{s}\right) $ such that $x\in V_{\alpha _{x}}\subseteq \tau _{\left(
i,j\right) }^{s}$-$\limfunc{int}\left( \tau _{\left( j,i\right) }^{s}\text{-}%
\limfunc{cl}\left( V_{\alpha _{x}}\right) \right) \subseteq U_{\alpha _{x}}$%
. Hence $X=\dbigcup\nolimits_{x\in X}V_{\alpha _{x}}$ and thus the family $%
\left\{ V_{\alpha _{x}}:x\in X\right\} $ forms a $\tau _{\left( i,j\right)
}^{s}$-open cover of $\left( X,\tau _{\left( 1,2\right) }^{s},\tau _{\left(
2,1\right) }^{s}\right) $.\ Since $\left( X,\tau _{\left( 1,2\right)
}^{s},\tau _{\left( 2,1\right) }^{s}\right) $ is $\left( \tau _{\left(
i,j\right) }^{s},\tau _{\left( j,i\right) }^{s}\right) $-nearly Lindel\"{o}%
f, there exists a countable subset of points $x_{1},\ldots ,x_{n},\ldots $
of $X$ such that $X=\dbigcup\nolimits_{n\in 
\mathbb{N}
}\tau _{\left( i,j\right) }^{s}$-$\limfunc{int}\left( \tau _{\left(
j,i\right) }^{s}\text{-}\limfunc{cl}\left( V_{\alpha _{x_{n}}}\right)
\right) \subseteq \dbigcup\nolimits_{n\in 
\mathbb{N}
}U_{\alpha _{x_{n}}}$. This shows that $\left( X,\tau _{\left( 1,2\right)
}^{s},\tau _{\left( 2,1\right) }^{s}\right) $ is $\tau _{\left( i,j\right)
}^{s}$-Lindel\"{o}f.
\end{proof}

\begin{corollary}
A bitopological space $\left( X,\tau _{\left( 1,2\right) }^{s},\tau _{\left(
2,1\right) }^{s}\right) $ is pairwise nearly Lindel\"{o}f if and only if $%
\left( X,\tau _{\left( 1,2\right) }^{s},\tau _{\left( 2,1\right)
}^{s}\right) $\ is Lindel\"{o}f.
\end{corollary}

From the Definition \ref{Def 2.6}, if the property $\mathcal{P}$ is not
bitopological property but it satisfies the condition $\left( X,\tau
_{1},\tau _{2}\right) $\ has the property $\mathcal{P}$\ if and only if $%
\left( X,\tau _{\left( 1,2\right) }^{s},\tau _{\left( 2,1\right)
}^{s}\right) $\ has the property $\mathcal{P}$, then the property $\mathcal{P%
}$ will be called pairwise semiregular invariant property. The following
theorem prove that $\left( i,j\right) $-nearly Lindel\"{o}f as well as
pairwise nearly Lindel\"{o}f property satisfying the pairwise semiregular
invariant property since $\left( i,j\right) $-nearly Lindel\"{o}f and
pairwise nearly Lindel\"{o}f are not $i$-topological property \cite%
{AdemBidi07-1} and bitopological property, respectively. This is because the 
$i$-continuity and $\left( i,j\right) $-$\delta $-continuity (resp.
continuity and $p$-$\delta $-continuity) are independent notions (see \cite%
{AdemBidi2010}).

\begin{theorem}
A bitopological space $\left( X,\tau _{1},\tau _{2}\right) $\ is $\left(
\tau _{i},\tau _{j}\right) $-nearly Lindel\"{o}f if and only if $\left(
X,\tau _{\left( 1,2\right) }^{s},\tau _{\left( 2,1\right) }^{s}\right) $\ is 
$\left( \tau _{\left( i,j\right) }^{s},\tau _{\left( j,i\right) }^{s}\right) 
$-nearly Lindel\"{o}f.
\end{theorem}

\begin{proof}
It is obvious by Theorem \ref{Th 4.1} and Proposition \ref{Prop 4.1}. \ 
\end{proof}

\begin{corollary}
A bitopological space $\left( X,\tau _{1},\tau _{2}\right) $\ is pairwise
nearly Lindel\"{o}f if and only if $\left( X,\tau _{\left( 1,2\right)
}^{s},\tau _{\left( 2,1\right) }^{s}\right) $\ is pairwise nearly Lindel\"{o}%
f\textit{.}
\end{corollary}

\begin{theorem}
\label{Th 4.3}\cite{SingAry71}\textbf{\ }If $\left( X,\tau _{1},\tau
_{2}\right) $\ is pairwise semiregular, then $\left( X,\tau _{1},\tau
_{2}\right) =\left( X,\tau _{\left( 1,2\right) }^{s},\tau _{\left(
2,1\right) }^{s}\right) $.\ 
\end{theorem}

The converse of Theorem \ref{Th 4.3} is also true by the definitions.

\begin{proposition}
Let $\left( X,\tau _{1},\tau _{2}\right) $\ be a pairwise semiregular space.
Then $\left( X,\tau _{1},\tau _{2}\right) $ is $\left( i,j\right) $-nearly
Lindel\"{o}f if and only if it is $i$-Lindel\"{o}f.
\end{proposition}

\begin{proof}
By Theorem \ref{Th 4.3}, $\left( X,\tau _{1},\tau _{2}\right) =\left( X,\tau
_{\left( 1,2\right) }^{s},\tau _{\left( 2,1\right) }^{s}\right) $. The
result follows immediately by Proposition \ref{Prop 4.1}.
\end{proof}

\begin{corollary}
Let $\left( X,\tau _{1},\tau _{2}\right) $\ be a pairwise semiregular space.
Then $\left( X,\tau _{1},\tau _{2}\right) $ is pairwise nearly Lindel\"{o}f
if and only if it is Lindel\"{o}f.
\end{corollary}

Unlike all types of pairwise Lindel\"{o}f properties, the $\left( i,j\right) 
$-almost Lindel\"{o}f, pairwise almost Lindel\"{o}f, $\left( i,j\right) $%
-weakly Lindel\"{o}f and pairwise weakly Lindel\"{o}f properties are
pairwise semiregular properties as we prove in the following theorems.

\begin{theorem}
\label{Th 4.4}A bitopological space $\left( X,\tau _{1},\tau _{2}\right) $\
is $\left( \tau _{i},\tau _{j}\right) $-almost Lindel\"{o}f if and only if $%
\left( X,\tau _{\left( 1,2\right) }^{s},\tau _{\left( 2,1\right)
}^{s}\right) $\ is $\left( \tau _{\left( i,j\right) }^{s},\tau _{\left(
j,i\right) }^{s}\right) $-almost Lindel\"{o}f.
\end{theorem}

\begin{proof}
Let $\left( X,\tau _{1},\tau _{2}\right) $ be a $\left( \tau _{i},\tau
_{j}\right) $-almost Lindel\"{o}f and let $\left\{ U_{\alpha }:\alpha \in
\Delta \right\} $ be a $\tau _{\left( i,j\right) }^{s}$-open cover of $%
\left( X,\tau _{\left( 1,2\right) }^{s},\tau _{\left( 2,1\right)
}^{s}\right) $. Since $\tau _{\left( i,j\right) }^{s}\subseteq \tau _{i}$, $%
\left\{ U_{\alpha }:\alpha \in \Delta \right\} $ is a $\tau _{i}$-open cover
of the $\left( \tau _{i},\tau _{j}\right) $-almost Lindel\"{o}f space $%
\left( X,\tau _{1},\tau _{2}\right) $. Then there is a countable subset $%
\left\{ \alpha _{n}:n\in 
\mathbb{N}
\right\} $ of $\Delta $ such that $X=\dbigcup\nolimits_{n\in 
\mathbb{N}
}\tau _{j}$-$\limfunc{cl}\left( U_{\alpha _{n}}\right) $. By Lemma \ref{Lem
2.1}, we have $X=\dbigcup\nolimits_{n\in 
\mathbb{N}
}\tau _{\left( j,i\right) }^{s}$-$\limfunc{cl}\left( U_{\alpha _{n}}\right) $%
, which implies $\left( X,\tau _{\left( 1,2\right) }^{s},\tau _{\left(
2,1\right) }^{s}\right) $ is $\left( \tau _{\left( i,j\right) }^{s},\tau
_{\left( j,i\right) }^{s}\right) $-almost Lindel\"{o}f.

Conversely suppose that $\left( X,\tau _{\left( 1,2\right) }^{s},\tau
_{\left( 2,1\right) }^{s}\right) $\textit{\ }is $\left( \tau _{\left(
i,j\right) }^{s},\tau _{\left( j,i\right) }^{s}\right) $-almost Lindel\"{o}f
and let $\left\{ V_{\alpha }:\alpha \in \Delta \right\} $ be a $\tau _{i}$%
-open cover of $\left( X,\tau _{1},\tau _{2}\right) $. Since $V_{\alpha
}\subseteq \tau _{i}$-$\limfunc{int}\left( \tau _{j}\text{-}\limfunc{cl}%
\left( V_{\alpha }\right) \right) $ and $\tau _{i}$-$\limfunc{int}\left(
\tau _{j}\text{-}\limfunc{cl}\left( V_{\alpha }\right) \right) \in $ $\tau
_{\left( i,j\right) }^{s}$, we have $\left\{ \tau _{i}\text{-}\limfunc{int}%
\left( \tau _{j}\text{-}\limfunc{cl}\left( V_{\alpha }\right) \right)
:\alpha \in \Delta \right\} $ is a $\tau _{\left( i,j\right) }^{s}$-open
cover of the $\left( \tau _{\left( i,j\right) }^{s},\tau _{\left( j,i\right)
}^{s}\right) $-almost Lindel\"{o}f space $\left( X,\tau _{\left( 1,2\right)
}^{s},\tau _{\left( 2,1\right) }^{s}\right) $. So there is a countable
subset $\left\{ \alpha _{n}:n\in 
\mathbb{N}
\right\} $ of $\Delta $ such that $X=\dbigcup\nolimits_{n\in 
\mathbb{N}
}$ $\tau _{\left( j,i\right) }^{s}$-$\limfunc{cl}\left( \tau _{i}\text{-}%
\limfunc{int}\left( \tau _{j}\text{-}\limfunc{cl}\left( V_{\alpha
_{n}}\right) \right) \right) $. By Lemma \ref{Lem 2.1}, we have $%
X=\dbigcup\nolimits_{n\in 
\mathbb{N}
}$ $\tau _{j}$-$\limfunc{cl}\left( \tau _{i}\text{-}\limfunc{int}\left( \tau
_{j}\text{-}\limfunc{cl}\left( V_{\alpha _{n}}\right) \right) \right)
\subseteq \dbigcup\nolimits_{n\in 
\mathbb{N}
}$ $\tau _{j}$-$\limfunc{cl}\left( V_{\alpha _{n}}\right) $. This implies
that $\left( X,\tau _{1},\tau _{2}\right) $ is $\left( \tau _{i},\tau
_{j}\right) $-almost Lindel\"{o}f.
\end{proof}

\begin{corollary}
\label{Cor 4.5}A bitopological space $\left( X,\tau _{1},\tau _{2}\right) $\
is pairwise almost Lindel\"{o}f if and only if $\left( X,\tau _{\left(
1,2\right) }^{s},\tau _{\left( 2,1\right) }^{s}\right) $\ is pairwise almost
Lindel\"{o}f.
\end{corollary}

Note that, the $\left( i,j\right) $-almost Lindel\"{o}f property and the
pairwise almost Lindel\"{o}f property are both bitopological properties (see 
\cite{BidiAdem2}). Utilizing this fact, Theorem \ref{Th 4.4} and Corollary %
\ref{Cor 4.5}, we easily obtain the following corollary.

\begin{corollary}
The $\left( i,j\right) $-almost Lindel\"{o}f property and the pairwise
almost Lindel\"{o}f property are both pairwise semiregular properties.
\end{corollary}

\begin{proposition}
Let $\left( X,\tau _{1},\tau _{2}\right) $\ be a $\left( \tau _{i},\tau
_{j}\right) $-almost regular space. Then $\left( X,\tau _{1},\tau
_{2}\right) $\ is $\left( \tau _{i},\tau _{j}\right) $-almost Lindel\"{o}f
if and only if $\left( X,\tau _{\left( 1,2\right) }^{s},\tau _{\left(
2,1\right) }^{s}\right) $\ is $\tau _{\left( i,j\right) }^{s}$-Lindel\"{o}f.
\end{proposition}

\begin{proof}
Let $\left( X,\tau _{1},\tau _{2}\right) $ be a $\left( \tau _{i},\tau
_{j}\right) $-almost Lindel\"{o}f and let $\left\{ U_{\alpha }:\alpha \in
\Delta \right\} $ be a $\tau _{\left( i,j\right) }^{s}$-open cover of $%
\left( X,\tau _{\left( 1,2\right) }^{s},\tau _{\left( 2,1\right)
}^{s}\right) $. For each $x\in X$, there exists $\alpha _{x}\in \Delta $
such that $x\in U_{\alpha _{x}}$ and since $U_{\alpha _{x}}\in $ $\tau
_{\left( i,j\right) }^{s}$, there exists a $\left( \tau _{i},\tau
_{j}\right) $-regular open set $V_{\alpha _{x}}$ in $\left( X,\tau _{1},\tau
_{2}\right) $ such that $x\in V_{\alpha _{x}}\subseteq U_{\alpha _{x}}$.
Since $\left( X,\tau _{1},\tau _{2}\right) $\textit{\ }is $\left( \tau
_{i},\tau _{j}\right) $-almost regular, there is a $\left( \tau _{i},\tau
_{j}\right) $-regular open set $C_{\alpha _{x}}$ in $\left( X,\tau _{1},\tau
_{2}\right) $ such that $x\in C_{\alpha _{x}}\subseteq \tau _{j}$-$\limfunc{%
cl}\left( C_{\alpha _{x}}\right) \subseteq V_{\alpha _{x}}$. Hence $%
X=\dbigcup\nolimits_{x\in X}C_{\alpha _{x}}$ and thus the family $\left\{
C_{\alpha _{x}}:x\in X\right\} $ forms a $\left( \tau _{i},\tau _{j}\right) $%
-regular open cover of $\left( X,\tau _{1},\tau _{2}\right) $. Since $\left(
X,\tau _{1},\tau _{2}\right) $ is $\left( \tau _{i},\tau _{j}\right) $%
-almost Lindel\"{o}f, there exists a countable subset of points $%
x_{1},\ldots ,x_{n},\ldots $ of $X$ such that $X=\dbigcup\nolimits_{n\in 
\mathbb{N}
}\tau _{j}$-$\limfunc{cl}\left( C_{\alpha _{x_{n}}}\right) \subseteq
\dbigcup\nolimits_{n\in 
\mathbb{N}
}V_{\alpha _{x_{n}}}\subseteq \dbigcup\nolimits_{n\in 
\mathbb{N}
}U_{\alpha _{x_{n}}}$. This shows that $\left( X,\tau _{\left( 1,2\right)
}^{s},\tau _{\left( 2,1\right) }^{s}\right) $ is $\tau _{\left( i,j\right)
}^{s}$-Lindel\"{o}f. Conversely, let $\left( X,\tau _{\left( 1,2\right)
}^{s},\tau _{\left( 2,1\right) }^{s}\right) $ be a $\tau _{\left( i,j\right)
}^{s}$-Lindel\"{o}f and let $\left\{ U_{\alpha }:\alpha \in \Delta \right\} $
be a $\tau _{i}$-open cover of $\left( X,\tau _{1},\tau _{2}\right) $. Since 
$U_{\alpha }\subseteq \tau _{i}$-$\limfunc{int}\left( \tau _{j}\text{-}%
\limfunc{cl}\left( U_{\alpha }\right) \right) $ and $\tau _{i}$-$\limfunc{int%
}\left( \tau _{j}\text{-}\limfunc{cl}\left( U_{\alpha }\right) \right) \in $ 
$\tau _{\left( i,j\right) }^{s}$, $\left\{ \tau _{i}\text{-}\limfunc{int}%
\left( \tau _{j}\text{-}\limfunc{cl}\left( U_{\alpha }\right) \right)
:\alpha \in \Delta \right\} $ is $\tau _{\left( i,j\right) }^{s}$-open cover
of the $\tau _{\left( i,j\right) }^{s}$-Lindel\"{o}f space $\left( X,\tau
_{\left( 1,2\right) }^{s},\tau _{\left( 2,1\right) }^{s}\right) $. Then
there exists a countable subset $\left\{ \alpha _{n}:n\in 
\mathbb{N}
\right\} $ of $\Delta $ such that $X=\dbigcup\nolimits_{n\in 
\mathbb{N}
}\tau _{i}$-$\limfunc{int}\left( \tau _{j}\text{-}\limfunc{cl}\left(
U_{\alpha _{n}}\right) \right) \subseteq \dbigcup\nolimits_{n\in 
\mathbb{N}
}\tau _{j}$-$\limfunc{cl}\left( U_{\alpha _{n}}\right) $. This implies that $%
\left( X,\tau _{1},\tau _{2}\right) $ is $\left( \tau _{i},\tau _{j}\right) $%
-almost Lindel\"{o}f.\bigskip
\end{proof}

\begin{corollary}
Let $\left( X,\tau _{1},\tau _{2}\right) $\ be a pairwise almost regular
space. Then $\left( X,\tau _{1},\tau _{2}\right) $\ is pairwise almost Lindel%
\"{o}f if and only if $\left( X,\tau _{\left( 1,2\right) }^{s},\tau _{\left(
2,1\right) }^{s}\right) $\ is Lindel\"{o}f.
\end{corollary}

\begin{proposition}
\label{Prop 4.4}Let $\left( X,\tau _{\left( 1,2\right) }^{s},\tau _{\left(
2,1\right) }^{s}\right) $\ be a $\left( \tau _{\left( j,i\right) }^{s},\tau
_{\left( i,j\right) }^{s}\right) $-extremally disconnected space. Then $%
\left( X,\tau _{\left( 1,2\right) }^{s},\tau _{\left( 2,1\right)
}^{s}\right) $ is $\left( \tau _{\left( i,j\right) }^{s},\tau _{\left(
j,i\right) }^{s}\right) $-almost Lindel\"{o}f if and only if $\left( X,\tau
_{\left( 1,2\right) }^{s},\tau _{\left( 2,1\right) }^{s}\right) $\ is $\tau
_{\left( i,j\right) }^{s}$-Lindel\"{o}f.
\end{proposition}

\begin{proof}
The sufficient condition is obvious by the definitions. So we need only to
prove necessary condition. Suppose that $\left\{ U_{\alpha }:\alpha \in
\Delta \right\} $ is a $\tau _{\left( i,j\right) }^{s}$-open cover of $%
\left( X,\tau _{\left( 1,2\right) }^{s},\tau _{\left( 2,1\right)
}^{s}\right) $. For each $x\in X$, there exists $\alpha _{x}\in \Delta $
such that $x\in U_{\alpha _{x}}$. Since $\left( X,\tau _{\left( 1,2\right)
}^{s},\tau _{\left( 2,1\right) }^{s}\right) $ is $\left( \tau _{\left(
i,j\right) }^{s},\tau _{\left( j,i\right) }^{s}\right) $-semiregular, there
exists a $\tau _{\left( i,j\right) }^{s}$-open set $V_{\alpha _{x}}$ in $%
\left( X,\tau _{\left( 1,2\right) }^{s},\tau _{\left( 2,1\right)
}^{s}\right) $ such that $x\in V_{\alpha _{x}}\subseteq \tau _{\left(
i,j\right) }^{s}$-$\limfunc{int}\left( \tau _{\left( j,i\right) }^{s}\text{-}%
\limfunc{cl}\left( V_{\alpha _{x}}\right) \right) \subseteq U_{\alpha _{x}}$%
. Hence $X=\dbigcup\nolimits_{x\in X}V_{\alpha _{x}}$ and thus the family $%
\left\{ V_{\alpha _{x}}:x\in X\right\} $ forms a $\tau _{\left( i,j\right)
}^{s}$-open cover of $\left( X,\tau _{\left( 1,2\right) }^{s},\tau _{\left(
2,1\right) }^{s}\right) $.\ Since $\left( X,\tau _{\left( 1,2\right)
}^{s},\tau _{\left( 2,1\right) }^{s}\right) $ is $\left( \tau _{\left(
i,j\right) }^{s},\tau _{\left( j,i\right) }^{s}\right) $-almost Lindel\"{o}f
and $\left( \tau _{\left( j,i\right) }^{s},\tau _{\left( i,j\right)
}^{s}\right) $-extremally disconnected, there exists a countable subset of
points $x_{1},\ldots ,x_{n},\ldots $ of $X$ such that $X=\dbigcup%
\nolimits_{n\in 
\mathbb{N}
}\tau _{\left( j,i\right) }^{s}$-$\limfunc{cl}\left( V_{\alpha
_{x_{n}}}\right) =\dbigcup\nolimits_{n\in 
\mathbb{N}
}\tau _{\left( i,j\right) }^{s}$-$\limfunc{int}\left( \tau _{\left(
j,i\right) }^{s}\text{-}\limfunc{cl}\left( V_{\alpha _{x_{n}}}\right)
\right) \subseteq \dbigcup\nolimits_{n\in 
\mathbb{N}
}U_{\alpha _{x_{n}}}$. This shows that $\left( X,\tau _{\left( 1,2\right)
}^{s},\tau _{\left( 2,1\right) }^{s}\right) $ is $\tau _{\left( i,j\right)
}^{s}$-Lindel\"{o}f.
\end{proof}

\begin{corollary}
Let $\left( X,\tau _{\left( 1,2\right) }^{s},\tau _{\left( 2,1\right)
}^{s}\right) $\ be a pairwise extremally disconnected space. Then $\left(
X,\tau _{\left( 1,2\right) }^{s},\tau _{\left( 2,1\right) }^{s}\right) $ is
pairwise almost Lindel\"{o}f if and only if $\left( X,\tau _{\left(
1,2\right) }^{s},\tau _{\left( 2,1\right) }^{s}\right) $\ is Lindel\"{o}f
\end{corollary}

\begin{proposition}
Let $\left( X,\tau _{1},\tau _{2}\right) $\ be a pairwise semiregular and $%
\left( j,i\right) $-extremally disconnected space. Then $\left( X,\tau
_{1},\tau _{2}\right) $ is $\left( i,j\right) $-almost Lindel\"{o}f if and
only if it is $i$-Lindel\"{o}f.
\end{proposition}

\begin{proof}
By Theorem \ref{Th 4.3}, $\left( X,\tau _{1},\tau _{2}\right) =\left( X,\tau
_{\left( 1,2\right) }^{s},\tau _{\left( 2,1\right) }^{s}\right) $. The
result follows immediately by Proposition \ref{Prop 4.4}.
\end{proof}

\begin{corollary}
Let $\left( X,\tau _{1},\tau _{2}\right) $\ be a pairwise semiregular and
pairwise extremally disconnected space. Then $\left( X,\tau _{1},\tau
_{2}\right) $ is pairwise almost Lindel\"{o}f if and only if it is Lindel%
\"{o}f.
\end{corollary}

\begin{theorem}
\label{Th 4.5}A bitopological space $\left( X,\tau _{1},\tau _{2}\right) $\
is $\left( \tau _{i},\tau _{j}\right) $-weakly Lindel\"{o}f if and only if $%
\left( X,\tau _{\left( 1,2\right) }^{s},\tau _{\left( 2,1\right)
}^{s}\right) $\ is $\left( \tau _{\left( i,j\right) }^{s},\tau _{\left(
j,i\right) }^{s}\right) $-weakly Lindel\"{o}f.
\end{theorem}

\begin{proof}
The proof is similar to the proof of Theorem \ref{Th 4.4} by using the fact
that%
\begin{eqnarray*}
\tau _{\left( j,i\right) }^{s}\text{-}\limfunc{cl}\left(
\dbigcup\nolimits_{n\in 
\mathbb{N}
}\tau _{i}\text{-}\limfunc{int}\left( \tau _{j}\text{-}\limfunc{cl}\left(
V_{\alpha _{n}}\right) \right) \right) &=&\tau _{j}\text{-}\limfunc{cl}%
\left( \dbigcup\nolimits_{n\in 
\mathbb{N}
}\tau _{i}\text{-}\limfunc{int}\left( \tau _{j}\text{-}\limfunc{cl}\left(
V_{\alpha _{n}}\right) \right) \right) \\
&\subseteq &\tau _{j}\text{-}\limfunc{cl}\left( \dbigcup\nolimits_{n\in 
\mathbb{N}
}\tau _{j}\text{-}\limfunc{cl}\left( V_{\alpha _{n}}\right) \right) \\
&\subseteq &\tau _{j}\text{-}\limfunc{cl}\left( \dbigcup\nolimits_{n\in 
\mathbb{N}
}V_{\alpha _{n}}\right) \text{.}
\end{eqnarray*}%
Thus we choose to omit the details.
\end{proof}

\begin{corollary}
\label{Cor 4.10}A bitopological space $\left( X,\tau _{1},\tau _{2}\right) $%
\ is pairwise weakly Lindel\"{o}f if and only if $\left( X,\tau _{\left(
1,2\right) }^{s},\tau _{\left( 2,1\right) }^{s}\right) $\ is pairwise weakly
Lindel\"{o}f.
\end{corollary}

Note that, the $\left( i,j\right) $-weakly Lindel\"{o}f property and the
pairwise weakly Lindel\"{o}f property are both bitopological properties (see 
\cite{BidiAdem2}). Utilizing this fact, Theorem \ref{Th 4.5} and Corollary %
\ref{Cor 4.10}, we easily obtain the following corollary.

\begin{corollary}
The $\left( i,j\right) $-weakly Lindel\"{o}f property and the pairwise
weakly Lindel\"{o}f property are both pairwise semiregular properties.
\end{corollary}

Recall that, a bitopological space $X$\ is called $\left( i,j\right) $-weak $%
P$-space \cite{BidiAdem1} if for each countable family $\left\{ U_{n}:n\in 
\mathbb{N}
\right\} $\ of $i$-open sets in $X$, we have $j$\textit{-}$\limfunc{cl}%
\left( \dbigcup\limits_{n\in 
\mathbb{N}
}U_{n}\right) =\dbigcup\limits_{n\in 
\mathbb{N}
}j$\textit{-}$\limfunc{cl}\left( U_{n}\right) $\textit{. }$X$\ is called
pairwise weak $P$-space if it is both $\left( 1,2\right) $-weak $P$-space
and $\left( 2,1\right) $-weak $P$-space.

\begin{proposition}
Let $\left( X,\tau _{1},\tau _{2}\right) $\ be a $\left( \tau _{i},\tau
_{j}\right) $-almost regular and $\left( \tau _{i},\tau _{j}\right) $-weak $%
P $-space. Then $\left( X,\tau _{1},\tau _{2}\right) $\ is $\left( \tau
_{i},\tau _{j}\right) $-weakly Lindel\"{o}f if and only if $\left( X,\tau
_{\left( 1,2\right) }^{s},\tau _{\left( 2,1\right) }^{s}\right) $\ is $\tau
_{\left( i,j\right) }^{s}$-Lindel\"{o}f.
\end{proposition}

\begin{proof}
Necessity: Let $\left\{ U_{\alpha }:\alpha \in \Delta \right\} $ be a $\tau
_{\left( i,j\right) }^{s}$-open cover of $\left( X,\tau _{\left( 1,2\right)
}^{s},\tau _{\left( 2,1\right) }^{s}\right) $. For each $x\in X$, there
exists $\alpha _{x}\in \Delta $ such that $x\in U_{\alpha _{x}}$ and since $%
U_{\alpha _{x}}\in $ $\tau _{\left( i,j\right) }^{s}$ for each $\alpha
_{x}\in \Delta $, there exists a $\left( \tau _{i},\tau _{j}\right) $%
-regular open set $V_{\alpha _{x}}$ in $\left( X,\tau _{1},\tau _{2}\right) $
such that $x\in V_{\alpha _{x}}\subseteq U_{\alpha _{x}}$. Since $\left(
X,\tau _{1},\tau _{2}\right) $\textit{\ }is $\left( \tau _{i},\tau
_{j}\right) $-almost regular, there is a $\left( \tau _{i},\tau _{j}\right) $%
-regular open set $C_{\alpha _{x}}$ in $\left( X,\tau _{1},\tau _{2}\right) $
such that $x\in C_{\alpha _{x}}\subseteq \tau _{j}$-$\limfunc{cl}\left(
C_{\alpha _{x}}\right) \subseteq V_{\alpha _{x}}$. Hence $%
X=\dbigcup\nolimits_{x\in X}C_{\alpha _{x}}$ and thus the family $\left\{
C_{\alpha _{x}}:x\in X\right\} $ forms a $\left( \tau _{i},\tau _{j}\right) $%
-regular open cover of $\left( X,\tau _{1},\tau _{2}\right) $. Since $\left(
X,\tau _{1},\tau _{2}\right) $ is $\left( \tau _{i},\tau _{j}\right) $%
-weakly Lindel\"{o}f and $\left( \tau _{i},\tau _{j}\right) $-weak $P$%
-space, there exists a countable subset of points $x_{1},\ldots
,x_{n},\ldots $ of $X$ such that $X=\tau _{j}$-$\limfunc{cl}\left(
\dbigcup\nolimits_{n\in 
\mathbb{N}
}C_{\alpha _{x_{n}}}\right) =\dbigcup\nolimits_{n\in 
\mathbb{N}
}\tau _{j}$-$\limfunc{cl}\left( C_{\alpha _{x_{n}}}\right) \subseteq
\dbigcup\nolimits_{n\in 
\mathbb{N}
}V_{\alpha _{x_{n}}}\subseteq \dbigcup\nolimits_{n\in 
\mathbb{N}
}U_{\alpha _{x_{n}}}$. This shows that $\left( X,\tau _{\left( 1,2\right)
}^{s},\tau _{\left( 2,1\right) }^{s}\right) $ is $\tau _{\left( i,j\right)
}^{s}$-Lindel\"{o}f.

Sufficiency: Let $\left\{ U_{\alpha }:\alpha \in \Delta \right\} $ be a $%
\tau _{i}$-open cover of $\left( X,\tau _{1},\tau _{2}\right) $. Since $%
U_{\alpha }\subseteq \tau _{i}$-$\limfunc{int}\left( \tau _{j}\text{-}%
\limfunc{cl}\left( U_{\alpha }\right) \right) $ and $\tau _{i}$-$\limfunc{int%
}\left( \tau _{j}\text{-}\limfunc{cl}\left( U_{\alpha }\right) \right) \in $ 
$\tau _{\left( i,j\right) }^{s}$, $\left\{ \tau _{i}\text{-}\limfunc{int}%
\left( \tau _{j}\text{-}\limfunc{cl}\left( U_{\alpha }\right) \right)
:\alpha \in \Delta \right\} $ is $\tau _{\left( i,j\right) }^{s}$-open cover
of the $\tau _{\left( i,j\right) }^{s}$-Lindel\"{o}f space$\left( X,\tau
_{\left( 1,2\right) }^{s},\tau _{\left( 2,1\right) }^{s}\right) $. Then
there exists a countable subset $\left\{ \alpha _{n}:n\in 
\mathbb{N}
\right\} $ of $\Delta $ such that $X=\dbigcup\nolimits_{n\in 
\mathbb{N}
}\tau _{i}$-$\limfunc{int}\left( \tau _{j}\text{-}\limfunc{cl}\left(
U_{\alpha _{n}}\right) \right) \subseteq \dbigcup\nolimits_{n\in 
\mathbb{N}
}\tau _{j}$-$\limfunc{cl}\left( U_{\alpha _{n}}\right) =\tau _{j}$-$\limfunc{%
cl}\left( \dbigcup\nolimits_{n\in 
\mathbb{N}
}U_{\alpha _{n}}\right) $. This implies that $\left( X,\tau _{1},\tau
_{2}\right) $ is $\left( \tau _{i},\tau _{j}\right) $-weakly Lindel\"{o}f.
\end{proof}

\begin{corollary}
Let $\left( X,\tau _{1},\tau _{2}\right) $\ be a pairwise almost regular and
pairwise weak $P$-space. Then $\left( X,\tau _{1},\tau _{2}\right) $\ is
pairwise weakly Lindel\"{o}f if and only if $\left( X,\tau _{\left(
1,2\right) }^{s},\tau _{\left( 2,1\right) }^{s}\right) $\ is Lindel\"{o}f.
\end{corollary}

\begin{proposition}
\label{Prop 4.7} Let $\left( X,\tau _{\left( 1,2\right) }^{s},\tau _{\left(
2,1\right) }^{s}\right) $\ be a $\left( \tau _{\left( j,i\right) }^{s},\tau
_{\left( i,j\right) }^{s}\right) $-extremally disconnected and $\left( \tau
_{\left( i,j\right) }^{s},\tau _{\left( j,i\right) }^{s}\right) $-weak $P$%
-space. Then $\left( X,\tau _{\left( 1,2\right) }^{s},\tau _{\left(
2,1\right) }^{s}\right) $ is $\left( \tau _{\left( i,j\right) }^{s},\tau
_{\left( j,i\right) }^{s}\right) $-weakly Lindel\"{o}f if and only if $%
\left( X,\tau _{\left( 1,2\right) }^{s},\tau _{\left( 2,1\right)
}^{s}\right) $\ is $\tau _{\left( i,j\right) }^{s}$-Lindel\"{o}f.
\end{proposition}

\begin{proof}
The sufficient condition is obvious by the definitions. So we need only to
prove necessary condition. Suppose that $\left\{ U_{\alpha }:\alpha \in
\Delta \right\} $ is a $\tau _{\left( i,j\right) }^{s}$-open cover of $%
\left( X,\tau _{\left( 1,2\right) }^{s},\tau _{\left( 2,1\right)
}^{s}\right) $. For each $x\in X$, there exists $\alpha _{x}\in \Delta $
such that $x\in U_{\alpha _{x}}$. Since $\left( X,\tau _{\left( 1,2\right)
}^{s},\tau _{\left( 2,1\right) }^{s}\right) $ is $\left( \tau _{\left(
i,j\right) }^{s},\tau _{\left( j,i\right) }^{s}\right) $-semiregular, there
exists a $\tau _{\left( i,j\right) }^{s}$-open set $V_{\alpha _{x}}$ in $%
\left( X,\tau _{\left( 1,2\right) }^{s},\tau _{\left( 2,1\right)
}^{s}\right) $ such that $x\in V_{\alpha _{x}}\subseteq \tau _{\left(
i,j\right) }^{s}$-$\limfunc{int}\left( \tau _{\left( j,i\right) }^{s}\text{-}%
\limfunc{cl}\left( V_{\alpha _{x}}\right) \right) \subseteq U_{\alpha _{x}}$%
. Hence $X=\dbigcup\nolimits_{x\in X}V_{\alpha _{x}}$ and thus the family $%
\left\{ V_{\alpha _{x}}:x\in X\right\} $ forms a $\tau _{\left( i,j\right)
}^{s}$-open cover of $\left( X,\tau _{\left( 1,2\right) }^{s},\tau _{\left(
2,1\right) }^{s}\right) $.\ Since $\left( X,\tau _{\left( 1,2\right)
}^{s},\tau _{\left( 2,1\right) }^{s}\right) $ is $\left( \tau _{\left(
i,j\right) }^{s},\tau _{\left( j,i\right) }^{s}\right) $-weakly Lindel\"{o}%
f, $\left( \tau _{\left( j,i\right) }^{s},\tau _{\left( i,j\right)
}^{s}\right) $-extremally disconnected and $\left( \tau _{\left( i,j\right)
}^{s},\tau _{\left( j,i\right) }^{s}\right) $-weak $P$-space, there exists a
countable subset of points $x_{1},\ldots ,x_{n},\ldots $ of $X$ such that $%
X=\tau _{\left( j,i\right) }^{s}$-$\limfunc{cl}\left(
\dbigcup\nolimits_{n\in 
\mathbb{N}
}V_{\alpha _{x_{n}}}\right) =\dbigcup\nolimits_{n\in 
\mathbb{N}
}\tau _{\left( j,i\right) }^{s}$-$\limfunc{cl}\left( V_{\alpha
_{x_{n}}}\right) =\dbigcup\nolimits_{n\in 
\mathbb{N}
}\tau _{\left( i,j\right) }^{s}$-$\limfunc{int}\left( \tau _{\left(
j,i\right) }^{s}\text{-}\limfunc{cl}\left( V_{\alpha _{x_{n}}}\right)
\right) \subseteq \dbigcup\nolimits_{n\in 
\mathbb{N}
}U_{\alpha _{x_{n}}}$. This shows that $\left( X,\tau _{\left( 1,2\right)
}^{s},\tau _{\left( 2,1\right) }^{s}\right) $ is $\tau _{\left( i,j\right)
}^{s}$-Lindel\"{o}f.
\end{proof}

\begin{corollary}
Let $\left( X,\tau _{\left( 1,2\right) }^{s},\tau _{\left( 2,1\right)
}^{s}\right) $\ be a pairwise extremally disconnected and pairwise weak $P$%
-space. Then $\left( X,\tau _{\left( 1,2\right) }^{s},\tau _{\left(
2,1\right) }^{s}\right) $ is pairwise weakly Lindel\"{o}f if and only if $%
\left( X,\tau _{\left( 1,2\right) }^{s},\tau _{\left( 2,1\right)
}^{s}\right) $\ is Lindel\"{o}f.
\end{corollary}

\begin{proposition}
Let $\left( X,\tau _{1},\tau _{2}\right) $\ be a pairwise semiregular, $%
\left( j,i\right) $-extremally disconnected and $\left( i,j\right) $-weak $P$%
-space. Then $\left( X,\tau _{1},\tau _{2}\right) $ is $\left( i,j\right) $%
-weakly Lindel\"{o}f if and only if it is $i$-Lindel\"{o}f.
\end{proposition}

\begin{proof}
By Theorem \ref{Th 4.3}, $\left( X,\tau _{1},\tau _{2}\right) =\left( X,\tau
_{\left( 1,2\right) }^{s},\tau _{\left( 2,1\right) }^{s}\right) $. The
result follows immediately by Proposition \ref{Prop 4.7}.
\end{proof}

\begin{corollary}
Let $\left( X,\tau _{1},\tau _{2}\right) $\ be a pairwise semiregular,
pairwise extremally disconnected and pairwise weak $P$-space. Then $\left(
X,\tau _{1},\tau _{2}\right) $ is pairwise weakly Lindel\"{o}f if and only
if it is Lindel\"{o}f.
\end{corollary}

\end{document}